\documentclass[11pt,leqno]{amsart}
\usepackage{blindtext}
\usepackage{fancyhdr}
\usepackage{ifpdf}
\usepackage{curves}
\usepackage{tikz}
\usetikzlibrary{arrows}
\input{arrowsnew}
\usetikzlibrary{decorations.markings}
\usepackage{subcaption}
\usepackage{placeins}
\usepackage{tabularx}
\usepackage{url}

\usepackage{amsfonts}
\usepackage{amsthm}
\usepackage{amsmath}
\usepackage{amscd}
\usepackage{amssymb}
\usepackage{graphicx}
\usepackage{rotating}
\usepackage[small,nohug,heads=LaTeX]{diagrams}
\diagramstyle[labelstyle=\scriptstyle]
\usepackage{ulem}
\usepackage{textcomp}
\usepackage{appendix}
\usepackage{chngcntr}
\usepackage{etoolbox}
\usepackage{lipsum}
\usepackage{tikz}
\usepackage{syntonly}
\usepackage[pagebackref]{hyperref}

\usepackage[utf8]{inputenc}
\usepackage[english, russian]{babel}

\usepackage{leftidx}
\usepackage{slashbox}

\usetikzlibrary{matrix,arrows,decorations.pathmorphing}

\AtBeginEnvironment{subappendices}{%
\chapter*{Appendix}
\addcontentsline{toc}{chapter}{Appendices}
\counterwithin{figure}{section}
\counterwithin{table}{section}
}

\def\CC {{\mathbb C}}     
\def\KK {{\mathbb K}}     
\def\NN {{\mathbb N}}     
\def\PP {{\mathbb P}}     
\def\XX {{\mathbb X}}     
\def\ZZ {{\mathbb Z}}     

\def\ring#1{\ifmmode \mathaccent'027 #1\else \rm\accent'027 #1\fi}

\newcommand{\rp}{{\mathbf p}}

\newcommand{\rL}{{\mathbf L}}

\def\ul  {\underline}

\def\wt  {\widetilde}

\def\mc {\mathcal}

\def\Hom {\mathrm{Hom}}

\def\coh {\mathrm{coh}}
\def\vect {\mathrm{vect}}
\def\mod {\mathrm{mod}}

\def \bd {\begin{diagram}}
\def \ed {\end{diagram}}
\def\be  {\begin{eqnarray}}
\def\ee  {\end{eqnarray}}
\def\ben {\begin{eqnarray*}}
\def\een {\end{eqnarray*}}

\def\bpr {\begin{proof}[Proof]}
\def\epr {\end{proof}}
\def\bsp {\begin{split}}
\def\esp {\end{split}}
\def\bprr {\begin{proof}[solution]}
\def\bpru {\begin{proof}[hint]}
\def\bpro {\begin{proof}[answer]}
\def\bcd {\begin{CD}}
\def\ecd {\end{CD}}

\newcommand{\abs}[1]{\left\vert#1\right\vert}

\newtheorem{theorem}{Theorem}[section]
\newtheorem{lemma}[theorem]{Lemma}
\newtheorem{fact}[theorem]{Fact}
\newtheorem{prop}[theorem]{Proposition}
\newtheorem{coro}[theorem]{Corollary}
\newtheorem{remark}[theorem]{Remark}
\newtheorem{df}[theorem]{Definition}

\newtheorem{conj}[theorem]{Conjecture}

\addcontentsline{toc}{chapter}{Appendices}
\begin{document}

\renewcommand{\contentsname}{Contents}
\renewcommand{\refname}{References}

\title[More finite sets coming from non-commutative counting]%
{More finite sets coming from non-commutative counting}
 
\author{George Dimitrov}
\address[Dimitrov]{Universit\"at Wien\\
	Oskar-Morgenstern-Platz 1, 1090 Wien\\
	\"Osterreich \\
}
\email{george.dimitrov@univie.ac.at}

\author{Ludmil Katzarkov}
\address[Katzarkov]{
 National Research University Higher School of Economics, Russian Federation, University of Miami  }
\email{lkatzarkov@gmail.com}

\renewcommand{\abstractname}{Abstract}

\renewcommand{\figurename}{Figure}

\begin{abstract} 

In our previous papers we introduced categorical invariants, which are, roughly speaking,   sets of triangulated  subcategories   in a given  triangulated category and their quotients.  Here  is extended the list of examples, where these sets are finite. Using results by Geigle, Lenzning, Meltzer, Hübner for weighted projective lines we show that for any two affine acyclic quivers $Q$, $Q'$  (i.e. quivers of extended Dynkin type) there are only finitely many full triangulated subctegories in $D^b(Rep_{\mathbb K}(Q))$, which are equivalent to  $D^b(Rep_{\mathbb K}(Q'))$, where ${\mathbb K}$ is an algebraically closed field.  Some of the numbers counting the elements in these finite sets are explicitly determined. 
 \end{abstract}

\maketitle
\setcounter{tocdepth}{1}
\tableofcontents

\section{Introduction}

 In \cite{DK4} we introduced some  categorical invariants via counting of fully faithful  functors.  In \cite{DK5}  we focused on examples and  extended  our studies   beyond counting.
 
  The  idea of non-commutative counting from  \cite{DK4} is unfolded here  in Section \ref{section main Definition}.  Namely,  for a triangulated  category $\mc T$, a subgroup $\Gamma \subset {\rm Aut}(\mc T)$ of  exact equivalences, and a choice of some restrictions $P$ on fully faithful functors   we  define  the set of subcategories of $\mc T$, which are equivalent to another chosen triangulated  category $\mc A$ via an equivalence  which satisfy  $P$, and modulo $\Gamma$.   The resulting  set of equivalence classes  of  subcategories in $\mc T$    we denote  by $C_{\mc A, P}^{\Gamma}(\mc T)$.   
 When the categories $\mc T$, $\mc A$ are $\KK$-linear for some field $\KK$, then choosing the property $P$ imposed on the functors to be $\KK$-linear, and restricting  $\Gamma$ to be subgroup of the group of $\KK$-linear auto-equivalences  $ {\rm Aut}_{\KK}(\mc T)$ results in  the set   $C_{\mc A, P}^{\Gamma}(\mc T)$, which we denote by  $C_{\mc A, \KK}^{\Gamma}(\mc T)$ (see Section  \ref{section main Definition}).
  
 Following Kontsevich-Rosenberg  \cite{KR} we denote $D^b(Rep_\KK(K(l+1)))$ by  $N\PP^l$  for $l\geq -1$, where $K(l+1)$ is the Kroncker quiver as explained in Section \ref{notations}. 
 For a $\KK$-linear category $\mc T$ and any $l\geq -1$ we denote $C_{N\PP^l, \KK}^\Gamma(\mc T)$ by  $C_{l}^{\Gamma}(\mc T)$ and refer to this set as the set of non-commutative curves, sometimes nc curves for short, of \textit{non-commutative genus} $l$ in $\mc T$ modulo  $\Gamma$. Note that   $N\PP^{-1}$ has homological dimension $0$ and one  should think of  the ``curves'' of non-commutative  genus $-1$ as a kind of  degenerate curves.
 
In \cite{DK5} we enriched our invariants  with  structures: the  inclusion of subcategories makes them partially ordered sets, and on the other hand  considering semi-orthogonal pairs of subcategories as edges amounts to directed graphs.  In addition to computing the non-commutative curve-counting invariants in $D^b(Q)$ for the affine quivers shown in \eqref{Q1}, for $ A_n$, and for  $D_4$ the paper  \cite{DK5}  contains  deriving of  formulas for counting of the subcategories of type $D^b(A_k)$ in $D^b(A_n)$,  whereas  for the two affine quivers and for $D_4$ we determined and counted all generated by an exceptional collection subcategories in $D^b(Q)$. It turns out that the problem for counting $D^b(A_k)$ in $D^b(A_n)$ has a geometric combinatorial parallel -  counting of maps between polygons.
 Estimating the numbers counting non-commutative curves in $D^b({\mathbb P}^2)$ modulo the group of auto-equivalences  we proved finiteness and that the exact determining of these numbers  leads to solution of Markov problem. Via homological mirror symmetry this gives a new approach to this problem. 
Regarding  the  structure of a partially ordered  set mentioned above we initiated in \cite{DK5} intersection theory of non-commutative curves focusing on the  case of non-commutative genus zero. 
 In the last section of \cite{DK5} we built an analogue of the classical curve complex, introduced by Harvey and Harrer and used in a spectacular way in Teichmueller theory and Thurston theory by Farb, Minsky, Mazur and many others.   The above-mentioned structure of a directed graph (and related simplicial complex)   is a categorical analogue of the curve complex.
The paper \cite{DK5}  contains pictures of  the graphs in many examples  and presents also an approach to Markov conjecture via counting of subgraphs in a graph associated with $D^b(\PP^2)$.

 Some of the examples, studied in \cite{DK5}, are as follows: Let $Q_1, Q_2$ be the quivers: \be \label{Q1} Q_1=  \begin{diagram}[1em]
 	&       &  2  &       &    \\
 	& \ruTo &    & \rdTo &       \\
 	1  & \rTo  &    &       &  3
 \end{diagram}   \ \qquad \qquad \ Q_2= \begin{diagram}[1em]
 2 &  \rTo  &  3    \\
 \uTo &        & \uTo     \\
 1   & \rTo  &    4\end{diagram}  \ee 

then for any other affine quiver $Q'$ and for $i=1,2$ we have $\abs{C^{\{\rm Id\}}_{D^b(Q')}(D^b(Q_i))}<\infty$,  all the non-zero cardinalities are computed in \cite{DK5}.

 In this note  we extend the  result about finiteness, $\abs{C^{\{\rm Id\}}_{D^b(Q')}(D^b(Q_i))}<\infty$, to the case when $Q_i$ is any affine acyclic quiver. In addition to that, 
for any two affine quivers $Q$, $Q'$ we determine whether   $C^{\{\rm Id\}}_{D^b(Q')}(D^b(Q))$ is empty or not. The main result is Theorem \ref{main result}. 
  
After that in Section \ref{remarks on the number} we prove that $\abs{C_1^{{\rm Aut}(D^b(Q))}(D^b(Q))}=1$ for any acyclic affine connected quiver $Q$ (Corollary \ref{with the full group}). In Corollary \ref{coro for numbers}  we deduce the number $\abs{C_1^{\{\rm Id\}}(D^b(Q))}$ for any acyclic affine connected quiver $Q$.

Important role for the proof of  Theorem \ref{main result}, Corollary \ref{with the full group}, and Corollary \ref{coro for numbers}  play   results by Geigle, Lenzning, Meltzer, Hübner for weighted projective lines  in  the papers  \cite{GL}, \cite{GL1}, \cite{Lenzing}, \cite{Meltzer1}, \cite{HuLe}, \cite{Hu}, \cite{LM}. 
 The proofs here rely often upon \cite[Corollary 1.9]{Orlov} and on \cite{Miyachi} as well.

 In the final Section \ref{final} we  make some comments and conjectures. We expect that with the help of the cited here   works of  Geigle, Lenzning, Meltzer, Hübner   finiteness  of  other  sets of the type $C^{\{\rm Id \}}_{D^b(\coh(\XX'))}(D^b(\coh(\XX)))$ can be proved, where $\XX$, $\XX'$ are weighted projective lines.  In  Section \ref{final} we note also  that whereas there are infinitely many non-commutative curves in $D^b(\PP^2)$ and infinitely many smooth projective curves in $\PP^2$ we have $C_{D^b(S)}^{\{\rm Id \}}(D^b(\PP^2)) = \emptyset$ for any smooth projective curve $S$ (see Lemma \ref{no commuta in P2} and the discussion before it).  The final conjecture \ref{final conjecture} discusses some threefolds $X$ and some smooth projective curves $S$ for which $C_{D^b(S)}^{\{\rm Id \}}(D^b(X)) \neq  \emptyset$.

   \textit{{\bf Acknowledgements:}}

   We are very grateful to  Helmut  Lenzing for his  kind and helpful E-mail response to our  question  regarding \eqref{hypothesis}.
   
   The first author  was  supported by  FWF Project P 29178-N35.

   The second author was
   supported by Simons research grant, NSF DMS 150908,
   DMS-1265230, Simons collaborative Grant - HMS, Simons
   investigator grant - HMS.
    The second author is partially
   supported by Laboratory of Mirror Symmetry NRU HSE, RF Government
   grant, ag. 14.641.31.0001

\section{Notations} \label{notations} 

We fix a universe and assume that the  set of objects and  the set of morphisms of any category we consider are elements of this universe. 

 The shift functor  in a triangulated category ${\mathcal T}$ is designated sometimes by $[1]$. We write $\langle  S \rangle  \subset \mc T$ for  the triangulated subcategory of $\mc T$ 
 generated by $S$, when $S \subset Ob(\mc T)$. 
  We write $\Hom^i(X,Y)$ for  $\Hom(X,Y[i])$.
  In this paper $\KK$ denotes an algebraically closed field.  
  If $\mc T$ is $\KK$-linear triangulated category  we write  $\hom^i(X,Y)$ for  $\dim_\KK(\Hom(X,Y[i]))$, where $X,Y\in \mc T$. 

A $\KK$-linear  triangulated category $\mc T$ is called  \textit{ proper} if $\sum_{i\in \ZZ} \hom^i(X,Y)<+\infty$ for any two objects $X,Y$ in $\mc T$.

An \textit{exceptional object} in a $\KK$-linear triangulated category  is an object $E\in \mc T$ satisfying $\Hom^i(E,E)=0$ for $i\neq 0$ and  $\Hom(E,E)=\KK $. We denote  by ${\mc T}_{exc}$ the set of all exceptional objects in $\mc T$,

An \textit{exceptional collection} is a sequence $\mc E = (E_0,E_1,\dots,E_n)\subset \mc T_{exc}$ satisfying $\hom^*(E_i,E_j)=0$ for $i>j$.    If  in addition we have $\langle \mc E \rangle = \mc T$, then $\mc E$ will be called a full exceptional collection. 

If an exceptional collection  $\mc E = (E_0,E_1,\dots,E_n)\subset \mc T_{exc}$ satisfies  $\hom^k(E_i,E_j)=0$ for any $i,j$ and for $k\neq 0$, then it is said to be \textit{strong exceptional collection}.

\textit{An abelian category $\mc A$ is said to be hereditary, if ${\rm Ext}^i(X,Y)=0$ for any  $X,Y \in \mc A$ and $i\geq 2$,  it is said to be of finite length, if it is Artinian and Noterian.}

By $Q$ we denote an acyclic quiver and   by  $D^b(Rep_\KK(Q))$, or just $D^b(Q)$, -  the derived category of the category of $\KK$-representations of $Q$,  this is a $\KK$-linear triangulated category.

For an integer $l\geq 0$ the $l$-Kronecker quiver  (the quiver with two vertices and  $l$ parallel  arrows)  will be denoted by  $K(l)$.

The number of elements of a finite set $X$ we denote by $\abs{X}$. If $X$ is infinite, we write  $\abs{X}=\infty$.

\begin{df}\label{SOD}  If $\mc T$ is a  triangulated category,  $\mc T_1$,  $\mc T_2$, $\dots$,  $\mc T_n$  are triangulated subcategories in it  satisfying the equalities   $\mc T = \langle \mc T_1, \mc T_2,\dots \mc T_n \rangle$ and  $\Hom(\mc T_j, \mc T_i)=0$ for $j>i$, then we say that   $\mc T = \langle \mc T_1, \mc T_2,\dots \mc T_n \rangle $ is a \ul{semi-orthogonal decomposition}. For short we write SOD instead of semi-orthogonal decomposition. 
\end{df}

\section{The sets $C_{\mc A, \KK}^{\Gamma}(\mc T)$, $C_{l}^{\Gamma}(\mc T)$} \label{section main Definition}
For  $\KK$-linear triangulated categories  $\mc A$, $\mc T$ and a subgroup  $\Gamma\subset {\rm Aut}_{\KK}(\mc T)$  of the group of $\KK$-linear auto-equivalences we defined in \cite[Definition 4.2]{DK5} a set $C_{\mc A, \KK}^{\Gamma}(\mc T)$. For short we write $C_{l}^{\Gamma}(\mc T)$ instead of $C_{N\PP^l, \KK}^{\Gamma}(\mc T)$. In this note we  consider only $\KK$-linear categories and we will omit writing $\KK$ in the notations $C_{\mc A, \KK}^{\Gamma}(\mc T)$. In particular,  here we need only the following:  

\begin{lemma} \cite[Lemma 4.5]{DK5} \label{bijection} Let $\mc A$, $\mc T$ be  $\KK$-linear triangulated categories.  Let  $\Gamma\subset {\rm Aut}_{\KK}(\mc T)$ be a subgroup of the group of $\KK$-linear auto-equivalences.
	There is a  bijection: 	\begin{gather}  \label{bijection in formula}
	C_{\mc A, \KK}^{\Gamma}(\mc T) \rightarrow\left  \{\mc B \subset \mc T:\begin{array}{l} \mc B \ \mbox{is a full tr. subcategory s.t.}  \\  \mbox{there exists a $\KK$-linear exact} \\ \mbox{ equivalence} \ \bd \mc A & \rTo^{}& \mc B \ed    \end{array}
	\right \}/\Gamma  \nonumber	\end{gather}
	the quotient  by  $\Gamma$ on the right hand side  means that  this is set of  equivalence classes of subcategories, where two subcategories  $\mc B_1$, $\mc B_2$ are in the same class iff there exists $[\beta] \in \Gamma$, s.t. $\beta(\mc B_1)=\mc B_2$. The meaning of ${\rm Im}(F)$ and $\beta(\mc B_1)$ is  explained in \cite[Lemma 3.6]{DK5}.	
\end{lemma}

\begin{prop} \cite[Proposition 5.5.]{DK5} \label{bijection 123}
	Let $\mc T$ be a triangulated category linear over $\KK$.  Let $\mc T$ admits an enhancement.  Then the bijcetion from Lemma \ref{bijection}  can be described as follows:
	\begin{gather}
\label{bijection 1} 	C_{l}^{\{\rm Id\}}(\mc T) \rightarrow\left  \{\mc A\subset \mc T:\begin{array}{l} \mc A \ \mbox{is a full tr. subcategory s.t.} \ \mc A =\langle E_1,E_2 \rangle \\  \mbox{for some strong exceptional pair} \ (E_1,E_2) \\ \mbox{s.t.} \ \hom(E_1,E_2)=l+1 \end{array} \right \} . 
	\end{gather}
\end{prop}

\section{Finiteness of $\abs{C^{\{\rm Id\}}_{D^b(Q')}(D^b(Q))}$ for affine acyclic  quivers $Q$, $Q'$}

\begin{df} \label{the categories indexed by p} By $\rp$ we will denote a sequence $\rp=(p_1,p_2,\dots, p_t)$, $t\geq 1$ of integers such that  $p_i \geq 1$ for $1\leq i\leq t$.

	 Following \cite[p. 319]{GL1} the  sequences  $\rp = (p_1, p_2,1)$, $ p_1 \geq 1, p_2\geq 1$,   $\rp = (2, 2, p)$ with $p\geq 2$, $\rp = (2, 3, 3)$, $\rp = (2, 3, 4)$, $\rp = (2, 3, 5)$ \underline{will be said to be of Dynkin type}.   
	
	For a sequence of Dynkin type $\rp$   we denote by  $\mc T(\rp)$ a triangulated category which is a derived category of a corresponding extended Dynkin quiver. More precisely, for $p \geq 1, q\geq 1$ we  will denote by $\mc T(p,q,1)$ the derived category of representations of a quiver, whose graph is an extended Dynkin diagram ot type  $A$ with $p+q$ vertices and orientation determined by the pair $p$, $q$ as shown in \cite[Figure 6]{Miyachi}. In \cite[p. 363]{Miyachi} is explained that for any acyclic  quiver whose underlying graph is an extended Dynkin diagram of type A we have $D^b(Q) \cong \mc T(p,q,1)$ for some $p \geq 1, q\geq 1$.
	
	When writing $D^b\left (\wt{D_n}\right )$ ($n\geq4$), $D^b\left (\wt{E_6} \right )$, $D^b\left (\wt{E_7} \right)$, $D^b\left (\wt{E_8} \right)$  we mean derived categories of representations of  quivers whose underlying graphs are   corresponding extended Dynkin diagrams and with any choice of orientation (changing the orientation gives equivalent categories \cite[Section 5]{Miyachi}). With this stipulation we denote	 
	\begin{gather} \mc T(2, 2, p)= D^b\left (\wt{D_{p+2}}\right ), \quad  \mc T(2, 3, 3)= D^b\left (\wt{E_6} \right ) ,\\  \mc T(2, 3, 4)= D^b\left (\wt{E_7} \right ), \quad   \mc T(2, 3, 5)= D^b\left (\wt{E_8} \right ).   \end{gather}
\end{df}
\begin{df} \label{dominates} Let $\rp$, $\rp'$ be two  sequences of length $3$.    We write $\rp'  \preccurlyeq \rp$ iff there exists  a permutation $\sigma$ of $\{1,2,3\}$ such that  ${p'_{\sigma(i)}} \leq {p_i}$ for $i=1,2,3$. 
\end{df}

\begin{remark} \label{permutation remark}
In the beginning of \cite[Section 9]{GL1} is explained  how to  view a sequence of natural numbers  and a normalized sequence of points in $\PP^1$ (i.e. the latter  has  the same length as the first sequence and starts with $\infty$, $0$, $1$) as a  weight function $\mathbf{w}: \PP^1\rightarrow \NN$.  For the case of  a sequence of length three   $\rp=(p_1,p_2,p_3)$ the  weight function $\mathbf{w}_{\rp}$  has values $\mathbf{w}_{\rp}(\infty)=p_1$, $\mathbf{w}_{\rp}(0)=p_2$, $\mathbf{w}_{\rp}(1)=p_2$ and $1$ elsewhere.    On  \cite[p. 321]{GL1} is defined what is meant by saying that a weight function $\mathbf{w}$ dominates another weight function   $\mathbf{v}$.  It can be shown that for two sequences $\rp$,  $\rp'$ of length  $3$   we have  $\rp' \preccurlyeq \rp$, as defined in Definition \ref{dominates}  iff the weight function $\mathbf{w}_{\rp}$ corresponding to  $\rp$ dominates the weight function $\mathbf{w}_{\rp'}$  corresponding to   $\rp'$.\footnote{to show this in one direction one uses that $PSL(2,\KK)$ acts $3$-transitively on $\PP^1$ ( $\KK$ is algebraically closed field). }
\end{remark}
\begin{remark} \label{reminder}
	In \cite{DK4} was shown that (see \cite[(7) and (8)]{DK4}):
	\begin{gather}
	\label{vanishings general}
	\mbox{If} \  l\geq 2 \ \mbox{and $Q$ is affine aciclic quiver}  \ \Rightarrow \abs{C_{l}^{\{\rm Id\}}(D^b(Q)) }=0. \\
	\label{vanishings general 2}  \mbox{If} \  l\geq 1 \ \mbox{and $Q$ is Dynkin  quiver}  \ \Rightarrow \abs{C_{l}^{\{\rm Id\}}(D^b(Q)) }=0.  \end{gather} 
	
\end{remark} 
We will use later the following 
\begin{lemma} \label{Lemma for C_2} Let $\rp$ be of Dynkin type.  From Remark \ref{reminder} we get  $\abs{C_{l}^{\{\rm Id\}}(\mc T(\rp)) }=0$ for $l\geq 2$.  We claim that furthermore,  $\abs{C_{1}^{\{\rm Id\}}(\mc T(\rp)) }\neq 0$.
	\end{lemma}
	\bpr \cite[Lemma 3.38]{Dimitrov} ensures  an exceptional pair $(E_1, E_2)\in \mc T(\rp)$, such that $\hom(E_1,E_2)=2$ and $\hom^i(E_1,E_2)=0$ for $i<0$. From \cite[Lemma 2.3.4]{Meltzer1} it follows that also  $\hom^i(E_1,E_2)=0$ for $i>0$. Now using  Proposition \ref{bijection 123} we deduce   $\langle E_1,E_2 \rangle\in C_{1}^{\{\rm Id\}}(\mc T(\rp))$.
	\epr
\begin{lemma} \label{Lemma for C_n} 
	Let $Q$ be an acyclic connected quiver, which is neither Dynkin nor extended  Dynkin.   Then    $\abs{C_{l}^{\{\rm Id\}}(D^b(Q)) }\neq 0$ for some $l\geq 2$.
\end{lemma}
\bpr We use \cite[Proposition 3.34]{DHKK} to get an exceptional pair $(E_1,E_2)\in D^b(Q)$, s.t. $\Hom(E_1,E_2)\geq 3$ and then  the same arguments as in Lemma \ref{Lemma for C_2} compete the proof.   \epr

\begin{coro} \label{infinitely many pairwise exc objects}
	For any acyclic connected quiver $Q$, which is not Dynkin, there exist infinitely many pairwise non-equivalent exceptional objects in $Rep_\KK(Q)$.
\end{coro}
\bpr Lemmas \ref{Lemma for C_n}, \ref{Lemma for C_2} ensure that  there exists a triangulated  subcategory $\mc T\subset  D^b(Q)$ equivalent to $D^b(K(l))$ for some $l\geq 2$. Since in   $D^b(K(l))$ there exists a sequence of exceptional objects $\{s_i\}_{i\in \ZZ}$ such that $s_i \not \cong s_j[k]$ for any $i\neq j$ and any $k\in \ZZ$ (see \cite[Section 5]{DK41}), such sequence exists in $D^b(Q)$ as well. Now the lemma follows since after shift each $s_i$ becomes an exceptional representation in $Rep_\KK(Q)$. \epr
\subsection{Domestic weighted projective lines} 

For any  sequence of  $\rp=(p_1,p_2,\dots, p_t)$, $t\geq 1$ of integers  and $p_i \geq 1$ in \cite{GL} were introduced  hom-finite  hereditary abelian  categories  (called \textit{the category of coherent sheaves on weighted projective lines} of weight $\rp$). In our considerations we will have $t= 3$ 
  and in this case the sequence  $\rp$ determines by the method in \cite{GL} a unique up to equivalence hereditary abelian category, which we  will denote by $ \coh(\XX(\rp))$, as in \cite{GL} we refer to  the objects of  $ \coh(\XX(\rp))$ as to sheaves on the weighted projective line $\XX(\rp)$. The full subcategory in  $ \coh(\XX(\rp))$ whose objects are sheaves of finite length, is denoted by  $ \coh_0(\XX(\rp))$. There is another full subcategory, whose elements are referred to as vector bundles on  $\XX(\rp)$ and this subcategory will be  denoted by   $ \vect(\XX(\rp))$. We list some facts for $\coh(\XX(\rp))$, which we need later: 

\begin{fact}{\rm (\cite[Proposition 2.4.]{GL})} \label{direct sum tor and torfree} Each object $X\in \coh(\XX(\rp))$ is a direct sum of an object in $\coh_0(\XX(\rp))$ and an object in  $\vect(\XX(\rp))$. 
\end{fact}

\begin{fact}{\rm (\cite[see p. 25]{Meltzer1})} \label{finitely many exc simples} There are only finitely many objects in $\coh_0(\XX(\rp))$ up to isomorphism, which are exceptional (since $\coh(\XX(\rp))$ is hereditary, an exceptional object in  $\coh(\XX(\rp))$ is the same as an object  $A\in \coh(\XX(\rp))$, s. t. $\Hom(X,X)=\KK$,  ${\rm Ext}^1(X,X)=0$).
\end{fact}

The  rank one abelian group  $\rL(\rp)$ on generators $\textbf{x}_1$, $\textbf{x}_2$, $\dots$, $\textbf{x}_t$ and relations $p_i \textbf{x}_i = p_j \textbf{x}_j$ for $1\leq i < j \leq t$ acts on $\coh(\XX(\rp))$ by exact auto-equivalences. For $\textbf{x}\in \rL(\rp)$  and $X\in \coh(\XX(\rp))$ the action of $\textbf{x}$ on $X$ is denoted by $X(\textbf{x})$. The element $p_1 \textbf{x}_1 = \dots = p_t \textbf{x}_t$ in $\rL(\rp)$ is denoted by $\textbf{c}$.

 One can define rank of the objects in $\coh(\XX(\rp))$, which is zero for the objects in $\coh_0(\XX(\rp))$  and positive integer on the objects in $\vect(\XX(\rp))$ (\cite[Corollary 1.8.2.]{GL}). 

\begin{fact}{\rm (\cite[ Proposition 2.1.]{GL})} \label{Pic} There is a specific rank one vector bundle, denoted by $\mc O \in \vect(\XX(\rp))$, such that the function $\rL(\rp) \ni x \mapsto \mc O (x)$ is a bijection between $\rL(\rp)$ and all rank-one vector bundles  in $\vect(\XX(\rp))$.
	\end{fact}

\begin{fact} (follows from  \cite[ Corollary 1.8.1, Proposition. 4.1]{GL}, see also \cite[p. 37, Lemma 2.3.4]{Meltzer1} ) \label{canonical exceptional collection} For any weight sequence   $\rp=(p_1,p_2,\dots, p_t)$  the sequence:
	\begin{gather} \label{canEXCcol} (\mc O, \mc O(\textbf{x}_1), \dots \mc O((p_1-1)\textbf{x}_1), \mc O(\textbf{x}_2), \dots \mc O((p_2-1)\textbf{x}_2), \dots, \mc O(\textbf{x}_t), \dots \mc O((p_t-1)\textbf{x}_t), \mc O(\textbf{c})  ) \end{gather}
	  is a full strong  exceptional collection on $D^b(\coh(\XX(\rp)))$. 
	\end{fact}
	
\begin{fact}(see \cite[Lemmas 3.1.2, 3.1.3]{Meltzer1}) \label{extendability of exceptional collections}  Each exceptional collection in $D^b(\coh(\XX(\rp)))$  of length equal to  ${\rm rank}(K_0(D^b(\coh(\XX(\rp)))))$ is full.  Each exceptional collection in $D^b(\coh(\XX(\rp)))$ can be extended to a full exceptional collection.

\end{fact}	

\begin{fact}(see \cite[Lemma 3.1.4]{Meltzer1}) \label{make an element unique} Let  $E_1,\dots, E_n$, $D_1,\dots, D_n$ be  two exceptional collections in $\coh(\XX(\rp))$  of length equal to  ${\rm rank}(K_0(D^b(\coh(\XX(\rp)))))$. If $E_i\cong D_i$ for $i\neq j$, then also $E_j\cong D_j$. 
\end{fact}		
	
\begin{fact} \label{the endomorphism algebra is canonical}
	Let $ E$ be the direct sum of the elements of  sequence \eqref{canEXCcol}.  Then the finite dimensional  $\KK$-algebra $\Hom(E,E)$ is the same as the so called canonical algebra of weight $\rp$ introduced in \cite{Ringel} and denoted by $\Lambda(\rp)$. The algebra  $\Lambda(\rp)$ is obtained via quiver with relations as explained in \cite[Section 4]{GL} or  \cite[p.  19,20]{Meltzer1}.
\end{fact}

\begin{fact}{\rm (\cite[see subsection 5.4.1]{GL})} \label{equivalences} For any weight  sequence $\rp$ we have a $\KK$-linear equivalence $ D^b(\coh(\XX(\rp)))\cong  D^b(\mod-\Lambda(\rp))$. Furthermore, if $\rp$ is of Dynikn type, in which case $\XX(\rp)$ is said to be of domestic type, we have also $\KK$-linear equivalence  $\mc T(\rp)\cong  D^b(\mod-\Lambda(\rp))$ ($\mc T(\rp)$ is explained in Definition \ref{the categories indexed by p}).\footnote{Note that for an acyclic affine quiver of type $A$ the category of its representations is equivalent to the category of representations of the dual quiver, so we we have these equivalences although we work with covariant representations of quivers in Definition \ref{the categories indexed by p}}
\end{fact}

\begin{fact} \label{equivalences for different ps} From \cite[Proposition 9.1]{GL1} and Remark \ref{permutation remark} it follows that if $\rp_1$ and $\rp_2$ are two weight sequences of length $3$ related by a permutation, then $\coh(\XX(\rp_1))$ is equivalent to  $\coh(\XX(\rp_2))$. 
\end{fact}

Using Facts \ref{equivalences for different ps} and \ref{equivalences} one shows that:

\begin{coro} \label{closedness property of sequivexes of Dynkin type}  Let $\rp$, $\rp'$ be two  sequences with length $3$.  Let $\rp$ be of Dynkin type. If $\rp'\preccurlyeq \rp$ (as defined in Definition \ref{dominates}), then up to permutation $\rp'$ is of Dynkin type and assuming that $\mathbf{q}$ is the permuted sequence, which is of Dynkin type, we have $ D^b(\coh(\XX(\rp')))\cong \mc T(\mathbf{q})$. 
\end{coro}

\begin{fact}{ \cite[PROPOSITION 1.1]{GL1}, here we use also that $\coh(\XX(\rp))$ is hereditary.} \label{perpedicular to a vector bundle} Let $E$ be an  exceptional sheaf on a weighted projective line $\XX(\rp)$. Then the full subcategory:
	\begin{gather} \label{the perp in coh} \{E\}^{\perp} = \{X \in \coh(\XX(\rp)): {\rm Ext}^i(E,X)=0 \ \forall i \} \end{gather}
	is closed under kernels, co-kernels and extensions, in particular it is an abelian subcategory in $\coh(\XX(\rp))$ and ${\rm Ext}_{\{E\}^{\perp}}^1(A,B) = {\rm Ext}_{\coh(\XX(\rp))}^1(A,B)$, $\Hom_{\{E\}^{\perp}}(A,B)=\Hom_{\coh(\XX(\rp))}(A,B)$ for any $A,B \in \{E\}^{\perp}$ (see the arguments in \cite[p. 5]{Lenzing}).

\end{fact}

\begin{remark} \label{reamrk for coproduct of categories}
	Let $\mc A$ be an abelian category. Let $\mc A_1$, $\mc A_2$ be two full subcategories, such that each object in $\mc A$ is a  direct sum of an object in $\mc A_1$ and an object in $\mc A_2$. Let $\Hom(X,Y)=\Hom(Y,X)=0$ for any $X\in \mc A_1$, and any $Y\in \mc A_2$. It follows that also ${\rm Ext}^1(X,Y)={\rm Ext}^1(Y,X)=0$ for any  $X\in \mc A_1$, and any $Y\in \mc A_2$. Furthermore, it follows that $\mc A_1$, $\mc A_2$ are Serre subcategories in $\mc A$ and  ${\rm Ext}_{\mc A_i}^1(A,B) = {\rm Ext}_{\mc A}^1(A,B)$, $\Hom_{\mc A_i}(A,B)=\Hom_{\mc A}(A,B)$ for any $A,B \in \mc A_i$ and any $i=1,2$.
\end{remark}

\begin{coro} \label{coro for per in coh} Let $E$ be an exceptional sheaf in  $\coh_0(\XX(\rp))$. Then the subcategory  $\{E\}^{\perp}$  in $\coh(\XX(\rp))$
	is coproduct of two subcategories $\mc A_1$, $\mc A_2$, such that $\mc A_1$ is equivalent to the category of representations of some quiver $Q$, whose  components are Dynkin quivers (i.e. of finite representation type), and $\mc A_2$ is equivalent to $\coh(\XX(\rp'))$ for some sequence $\rp'\prec\rp$ (strictly smaller for the defined in Definition \ref{dominates} partial order). Furthermore $\mc A_1\subset \coh_0(\XX(\rp))$. 
\end{coro}
\bpr  From  \cite[page 2550]{LM} (see also \cite[THEOREM 10.1 and the remark that follows]{GL1}) we get that  the subcategory: $\{E\}^{\perp}$ from \eqref{the perp in coh} in $\coh(\XX(\rp))$
is coproduct  of two subcategories $\mc A_1$, $\mc A_2$, such that $\mc A_1$ is equivalent to $\mod(H)$  for some hereditary finite dimensional $\KK$-algebra, and $\mc A_2$ is equivalent to $\coh(\XX(\rp'))$ for some sequence $\rp'\prec \rp$ (strictly smaller for the defined in Definition \ref{dominates} partial order). Furthermore $\mc A_1\subset \coh_0(\XX(\rp))$. 

It remains  to show that the category $\mod(H)$ is equivalent to the category of representations of some quiver $Q$ as specified. Since $H$ is finite dimensional hereditary it follows that $\mod(H)$ is equivalent to the category of representations of some acyclic quiver $Q$ (see e.g. \cite[p. 80]{Ringel}).   From Remark \ref{reamrk for coproduct of categories} and Fact \ref{perpedicular to a vector bundle} we see that any exceptional object in $\mc A_1$  is an exceptional object in $\coh(\XX(\rp))$ and since $\mc A_1$ is a full subcategory in $\coh(\XX(\rp))$ pairwise inequivalent exceptional objects in $\mc A_1$  are pairwise inequivalent in $\coh(\XX(\rp))$. Now from  Fact \ref{finitely many exc simples} and  $\mc A_1\subset \coh_0(\XX(\rp))$ it follows that there are finitely many exceptional objects in $Rep_\KK(Q)$. Now the corollary follows from Corollary \ref{infinitely many pairwise exc objects}.
\epr
\begin{coro}\label{perp to fl} Let $E$ be an exceptional sheaf in  $\coh_0(\XX(\rp))$. Then   we have a SOD $D^b(\coh(\XX(\rp))=\langle \mc T_1,\mc T_2,\langle  E \rangle \rangle$, where $\mc T_1 \cong D^b(Q)$ for some acyclic quiver $Q$ with components Dynkin quivers and  $\mc T_2 \cong D^b(\coh(\XX(\rp')))$ for some  $\rp'\prec \rp$. Furthermore $\Hom(\mc T_1, \mc T_2)=\Hom(\mc T_2, \mc T_1)=0$. 
	
	If $\rp$ is of Dynkin type, then    $\rp'$ can be chosen to be of Dynkin type.
\end{coro}	
\bpr  Provided that everything else is proved, the last sentence follows from Corollary \ref{closedness property of sequivexes of Dynkin type} and Fact \ref{equivalences}.

 Let $\mc A_1$, $\mc A_2$ be the subcategories in $\coh(\XX(\rp))$ from  Corollary \ref{coro for per in coh}. Using Remark \ref{reamrk for coproduct of categories} and  Fact \ref{perpedicular to a vector bundle} we see that exceptional collections in $\mc A_i$ give rise to exceptional collections in $D^b(\coh(\XX(\rp)))$ right orthogonal to $\langle E \rangle$. Furthermore from the same statements it follows that  for any object $X$ in $\mc A_1$ and any object $Y$ in $\mc A_2$ the pair $(X,Y)$ is orthogonal in $D^b(\coh(\XX(\rp)))$.

Let $A_1,\dots, A_n$ be the full  exceptional sequence in $D^b(\mc A_1)$ consisting of the projective  representations in $Rep_\KK(Q)$, this is a strong exceptional collection. Let $B_1,B_2,\dots,B_m$ be the full  exceptional collection in $D^b(\mc A_2)$ coming from the equivalence of $\mc A_2$  to $\coh(\XX(\rp'))$ and from Fact \ref{canonical exceptional collection}. We claim that  $A_1,\dots, A_n,B_1,\dots,B_m, E$ is a full exceptional collection in $D^b(\coh(\XX(\rp)))$. We already explained that   $\mc E=(A_1,\dots, A_n,B_1,\dots,B_m, E)$ is an exceptional collection in $D^b(\coh(\XX(\rp)))$ such that   $(A_i,B_j)$ are orthogonal pairs for any $i,j$.  If  $\mc E$ is not full then from Fact \ref{extendability of exceptional collections} there exists an exceptional object $X\in \langle \mc E \rangle^\perp$. Since $X$ is indecomposable we can assume that $X\in \coh(\XX(\rp))$ and then from  Corollary \ref{coro for per in coh} we have $X\in \mc A_i$ for $i=1$ or $i=2$. We assume $i=1$ (otherwise the arguments are the same). However using again Remark \ref{reamrk for coproduct of categories},  Fact \ref{perpedicular to a vector bundle}, and  $X\in \langle \mc E \rangle^\perp$  we deduce that ${\rm Ext}^j_{\mc A_1}(A_l ,X)=0$ for all $l$ and all $j$, which contradicts the fact that $(A_1,A_2,\dots,A_n)$ is a full exceptional collection in  $D^b(\mc A_1)$. So we proved that $\mc E$ is a full exceptional collection. Let us denote $\mc T_1 = \langle A_1,\dots, A_n \rangle$, $\mc T_2 = \langle B_1,\dots, B_m \rangle$, then from the already explained we see that  $D^b(\coh(\XX(\rp))=\langle \mc T_1,\mc T_2, E  \rangle$ is a SOD and $\Hom(\mc T_1, \mc T_2)=\Hom(\mc T_2, \mc T_1)=0$. Now both $(A_1,\dots, A_n)$, $(B_1,\dots, B_m)$ are full strong exceptional collections in $\mc T_1$, $\mc T_2$, respectively. Since $\mc A_1$ is equivalent to $Rep_\KK(Q)$ we have a full strong exceptional collection $(A_1',\dots, A_n')$ in $D^b(Q)$ such that the endomorphism algebra of $A_1'\oplus \dots \oplus A_n' $  is isomorphic to the endomorphism   algebra of $A_1\oplus \dots \oplus A_n$,  from  \cite[Corollary 1.9]{Orlov} it follows that $\mc T_1$ is equivalent to $D^b(Q)$.  By similar arguments and from the way we constructed the strong exceptional collection  $B_1,\dots,B_m$ it follows that $\mc T_2$ is equivalent to $D^b(\coh(\XX(\rp')))$, we take into account Facts \ref{the endomorphism algebra is canonical}, \ref{equivalences}.
\epr
We pass to studying the right orthogonal $\langle E \rangle^{\perp}$ in $D^b(\coh(\XX(\rp)))$ for an exceptional vector bundle $E$ on $\XX(\rp)$.

\begin{coro} \label{perpendicular to an exceptional vb} Let $\rp$ be of Dynkin type. Let $E$ be an exceptional vector bundle on  $\XX(\rp)$. Then  the subcategory  $\langle E \rangle^{\perp}$ in $D^b(\coh(\XX(\rp))$  is equivalent to $ D^b(Q)$ for some acyclic quiver $Q$ with $n-1$ vertices, where $n= {\rm rank}\left (K_0(\coh(\XX(\rp)))\right )$, with components Dynkin quivers. Furthermore, if $E$ is a line bundle, then $Q$ is Dynkin of type $A$, $D$, $E$, when $\rp$ is of type $A$, $D$, $E$, respectively (see Definition \ref{the categories indexed by p}). 
\end{coro}	
\bpr In this case the subcategory \eqref{the perp in coh} in $\coh(\XX(\rp))$ is equivalent to  $\mod(H)$ for some finite dimensional  hereditary $\KK$-algebra with $n-1$ simple objects  \cite[Theorem 2.4.3]{Meltzer1} (see also   the beginning of \cite[Section 2]{Meltzer1}, where is specified that by an algebra they always mean a finite dimensional algebra).  \cite[Theorem 2.4.3]{Meltzer1} is taken from \cite{HuLe}, it is  stated in \cite{Hu} as well.

Since $H$ is finite dimensional hereditary it follows that $\mod(H)$ is equivalent to the category of representations of some acyclic quiver $Q$ (see e.g. \cite[p. 80]{Ringel}).   From  Fact \ref{perpedicular to a vector bundle} we see that any exceptional collection  in $\{ E \}^{\perp}$ gives rise to an exceptional collection   in $\langle E \rangle^{\perp} \subset D^b(\coh(\XX(\rp)))$. Let $(E_1, E_2,\dots, E_{n-1})$ be the exceptional collection of  $\langle E \rangle^{\perp}$ coming from the projective representations in $Rep_\KK(Q)$, then by the already explained $(E_1, E_2,\dots, E_{n-1},E)$ is an exceptional collection in $D^b(\coh(\XX(\rp)))$ and from Fact \ref{extendability of exceptional collections} it is full, hence $\langle E \rangle^{\perp}=\langle E_1, E_2,\dots, E_{n-1} \rangle$. 
 Now as in the proof of the equivalence $\mc T_1 \cong D^b(Q)$ in  Corollary \ref{perp to fl} on proves that $\langle E \rangle^{\perp}\cong D^b(Q)$. Applying the explanation  after \cite[Theorem 2.4.3]{Meltzer1} to each of the possible $\rp$ one sees that if $E$ is a line bundle, then the quiver $Q$ is connected Dynkin of the corresponding type. 
 
 Consider  exceptional vector bundle $E$, which is not line bundle. Assume that  some component of $Q$ is not Dynkin.  From Lemmas \ref{Lemma for C_2}, \ref{Lemma for C_n}  $C_l\left ( \langle E \rangle^{\perp} \right )\neq \emptyset $ for some $l\geq 1$. Let $\langle E_1, E_2 \rangle \in C_l\left ( \langle E \rangle^{\perp} \right ) $, where $( E_1, E_2)$ is a strong  exceptional pair with $\hom(E_1,E_2)=l+1\geq 2$ due to Proposition \ref{bijection 123}. After shifting we obtain an exceptional pair $(X,Y)$ in $\{ E \}^{\perp}\subset \coh(\XX(\rp))$ with either $\hom(X,Y)=l+1$ or $\hom^1(X,Y)=l+1$.  Using again Fact \ref{extendability of exceptional collections} we obtain a complete exceptional collection $(X,Y,E_1,E_2,\dots,E_{n-3}, E)$ in $\coh(\XX(\rp))$ (by complete here we mean of length $n$). From \cite[Theorem 2.]{LM} and since $E$ is a vector bundle,  the smallest  exact (by exact subcategory we mean closed under kernels, co-kernels and finite  direct sums) and extension-closed  subcategory in $\coh(\XX(\rp))$ containing $(X,Y)$ is equivalent to a module category over a finite dimensional hereditary
 algebra, hence it  is equivalent to the category of representations of some acyclic quiver $Q$ (see e.g. \cite[p. 80]{Ringel}).
  Denote this subcategory by $\mc A$. Since it is  exact  and extension-closed, we have  for  any $A,B \in \mc A$:
 \begin{gather} \label{Homs and exts the same}  \Hom_{\coh(\XX(\rp))}(A,B) = \Hom_{\mc A}(A,B) \qquad {\rm Ext }^1_{\coh(\XX(\rp))}(A,B) = {\rm Ext }^1_{\mc A}(A,B)  \\ 
 \label{Hom 2 or ext 2} \Rightarrow \qquad  \Hom_{\mc A}(X,Y) \geq  2 \ \mbox{or} \  {\rm Ext }^1_{\mc A}(X,Y)\geq 2
 \end{gather}
 From the arguments in \cite[Proposition 4.2]{LM} we see that ${\rm rank}\left (K_0(\mc A )\right ) = 2$ and hence $\mc A $ is equivalent to $Rep_\KK(K(r))$ for some $r\geq 0$. Since all exceptional pairs in $Rep_\KK(K(r))$ have $r$-dimensional Hom or Ext between their elements,\footnote{see for example \cite[Figures 10,11,12 ]{DK5} or Remark \ref{remark for exceptianal pairs} below.}  from \eqref{Hom 2 or ext 2} it follows that $r=l+1\geq 2$. However this contradicts \cite[Proposition 4.3]{LM}. 
 
\epr

\begin{fact} (see \cite[THEOREM 9.8]{GL1} and Remark \ref{permutation remark}). \label{p' leq p} Let $\rp$, $\rp'$ be of Dynkin type and let $\rp'\preccurlyeq \rp$. Then $\coh(\XX(\rp'))$ is equivalent to a full exact subcategory of $\coh(\XX(\rp))$, which is closed under extensions.
\end{fact}

\begin{coro} \label{part of an important Corollary}  Let $\rp$, $\rp'$ be of Dynkin type and let $\rp'\preccurlyeq \rp$, then $C^{\{\rm Id\}}_{\mc T(\rp')}(\mc T(\rp)) \neq \emptyset$.
\end{coro}	
\bpr  Due to Fact \ref{equivalences} it is enough to show that $C^{\{\rm Id\}}_{D^b(\coh(\XX(\rp')))}(D^b(\coh(\XX(\rp)))) \neq \emptyset$. Using Fact \ref{p' leq p} we get  an exact and extension closed subcategory $\mc A \subset \coh(\XX(\rp))$, which is equivalent to $\coh(\XX(\rp'))$.  Since $\mc A$ is  exact  and extension-closed, we have  \eqref{Homs and exts the same} for  any $A,B \in \mc A$ again. This in turn ensures that the exceptional sequence on $\coh(\XX(\rp'))$ explained in Fact \ref{canonical exceptional collection} gives rise to a strong exceptional sequence $(E_1,\dots, E_n)$ in $D^b(\coh(\XX(\rp)))$, whose endomorphism algebra $\Hom\left (\oplus_{i=1}^n E_i, \oplus_{i=1}^n E_i\right )$ is $\Lambda(\rp')$ (see Fact \ref{the endomorphism algebra is canonical} as well).  Now \cite[Corollary 1.9]{Orlov} ensures that $\langle E_1,\dots, E_n\rangle \cong D^b(\mod-\Lambda(\rp')) $. Looking  again Fact \ref{equivalences} we  get $\langle E_1,\dots, E_n\rangle \cong D^b(\coh(\XX(\rp'))) $ and the corollary follows from Lemma \ref{bijection}. 
\epr

\begin{theorem} \label{main result} Let $\rp'$, $\rp$ be any two weight sequences of Dynkin type as defined in Definition \ref{the categories indexed by p}. Then $C^{\{\rm Id\}}_{\mc T(\rp')}(\mc T(\rp)) \neq \emptyset$ iff $\rp'\preccurlyeq \rp$ as defined in Definition \ref{dominates}. Furthermore  $\abs{C^{\{\rm Id\}}_{\mc T(\rp')}(\mc T(\rp))}<\infty$.
	
	It follows that for any two affine acyclic quivers $Q'$, $Q$ (i.e. quivers of extended Dynkin type) holds $\abs{C^{\{\rm Id\}}_{D^b(Q')}(D^b(Q))}<\infty$.
	\end{theorem}
	\bpr 
	The last sentence follows from Fact \ref{equivalences} and the comments in Definition \ref{the categories indexed by p}. 
	
	Taking  into account Corollary \ref{part of an important Corollary},  it is enough to prove that
	
	 \begin{gather} \label{finiteness formula} C^{\{\rm Id\}}_{\mc T(\rp')}(\mc T(\rp)) \neq \emptyset \ \  \mbox{\textit{implies}} \ \   \rp'\preccurlyeq \rp \ \ \mbox{\textit{ and}}  \ \ \abs{C^{\{\rm Id\}}_{\mc T(\rp')}(\mc T(\rp))}<\infty. \end{gather}

		Due to Fact \ref{equivalences} we can work with  $C^{\{\rm Id\}}_{D^b(\coh(\XX(\rp')))}(D^b(\coh(\XX(\rp))))$ instead of $C^{\{\rm Id\}}_{\mc T(\rp')}(\mc T(\rp)) $.  We  note that:
		
		\begin{gather}  \mbox{\textit{	Both $D^b(\coh(\XX(\rp')))$, $D^b(\coh(\XX(\rp)))$ have full exceptional collections, }} \nonumber  \\ \label{nice properties} \mbox{\textit{ each exceptional collection of maximal length is full, and each exceptional collection}}\\ \mbox{  \textit{can be extended to a full exceptional collection  (Facts \ref{canonical exceptional collection}, \ref{extendability of exceptional collections}}).} \nonumber\end{gather}

		\begin{gather} \mbox{\textit{Due to \eqref{nice properties} each $\mc A \in C^{\{\rm Id\}}_{D^b(\coh(\XX(\rp')))}(D^b(\coh(\XX(\rp))))$ is a subcategory} } \nonumber \\  \mbox{ \textit{generated by an exceptional collection and \cite[Theorem 3.2]{B} says  that}}\label{nice properties2}  \\ \mbox{  \textit{ $\mc A$ is admissible, which in turn implies that $\mc A$ is closed under direct summands.}} \nonumber \end{gather}

	Let us denote $\abs{\rp}=p_1+p_2+p_3$ for any sequence $\rp=(p_1,p_2,p_3)$.  We will make induction on $\abs{\rp}$.  We note first that from the very Definitions  \ref{the categories indexed by p}, \ref{dominates} we have	
	\begin{gather}\label{help for finiteness 0}  {\rm rank} \left ( \mc T(\rp) \right )= \abs{\rp}-1, \qquad   \rp'\prec \rp \Rightarrow \abs{\rp'} <  \abs{\rp'}.
		\end{gather}
	
		Recalling Definition \ref{the categories indexed by p} we see that the minimal  $\abs{\rp}$ is $3$ and there is a  unique $\rp$ where this is attained, it is  $\rp =(1,1,1)$. 
		
		If $\rp =(1,1,1)$ and  $C^{\{\rm Id\}}_{\mc T(\rp')}(\mc T(\rp)) \neq \emptyset$, then \eqref{help for finiteness 0}, \eqref{nice properties}, \eqref{nice properties2} imply  that $\abs{\rp'} \leq  \abs{\rp}=3$, and by the minimality of $\abs{\rp}$ and the  uniqueness of  $\rp$ we have $\rp' =  \rp$ and $\abs{C^{\{\rm Id\}}_{\mc T(\rp')}(\mc T(\rp))}=1$.

	The induction step: Let for some $N\geq 3$ we have that for any 
		$ \rp $   of Dynkin type and    with $\abs{\rp}\leq N $ \eqref{finiteness formula} holds.
	 Assume also that $C^{\{\rm Id\}}_{\mc T(\rp')}(\mc T(\rp)) \neq \emptyset$,  and that  $\abs{\rp}=N+1$.  We will show  that $\rp'\preccurlyeq \rp$ and that  $\abs{C^{\{\rm Id\}}_{\mc T(\rp')}(\mc T(\rp))}=\abs{C^{\{\rm Id\}}_{D^b(\coh(\XX(\rp')))}(D^b(\coh(\XX(\rp))))}<\infty$  and the theorem  will be proved. 
		
		If ${\rm rank}\left (K_0(D^b(\coh(\XX(\rp'))))\right)={\rm rank}\left (K_0(D^b(\coh(\XX(\rp))))\right)$, then from \eqref{nice properties}  it follows \\ $D^b(\coh(\XX(\rp')))\cong D^b(\coh(\XX(\rp))$, $\mc T(\rp)\cong \mc T(\rp')$ and from \cite[p. 56]{Ha} it follows that $\rp=\rp'$,  hence $\abs{C^{\{\rm Id\}}_{\mc T(\rp')}(\mc T(\rp))}=1$. 
		
		So, let ${\rm rank}\left (K_0(D^b(\coh(\XX(\rp'))))\right)<{\rm rank}\left (K_0(D^b(\coh(\XX(\rp))))\right)$. In this case \eqref{nice properties} ensures that for each $\mc A \in C^{\{\rm Id\}}_{D^b(\coh(\XX(\rp')))}(D^b(\coh(\XX(\rp))))$ there exists an exceptional sheaf $E\in \coh(\XX(\rp))$, s.t. $\mc A \subset \langle E \rangle^{\perp}$,  however Corollary \ref{perpendicular to an exceptional vb} shows that $E$ cannot be a vector bundle. Indeed if $E$ is a vector bundle, then from Corollary \ref{perpendicular to an exceptional vb}  and  $\mc A \subset \langle E \rangle^{\perp}$ would follow that for some acyclic quiver $Q$ with  components Dynkin quivers we have a subcategory $\mc A \subset D^b(Q)$ which is equivalent to $D^b(\coh(\XX(\rp')))\cong \mc T(\rp')$, which is impossible, since in $D^b(Q)$ there are only finitely many exceptional objects up to shift, whereas in  $\mc T(\rp')$, there are infinitely many indecomposable objects up to shift (we use \eqref{nice properties2} as well here).
				So far, we deduce that $\mc A \subset \langle E \rangle^{\perp}$ for some exceptional sheaf of finite length $E\in \coh_0(\XX(\rp))$, recall also Fact \ref{direct sum tor and torfree}.  From Fact \ref{finitely many exc simples} we know that the set of such objects, say $\mc E$, is finite. For each $ E \in \mc E$ from Corollary \ref{perp to fl}    we have $\langle E \rangle^{\perp} = \langle \mc T_1^{E}, \mc T_2^{E}\rangle  $ with $\KK$-linear equivalences $\mc T_1^E \cong D^b(Q_E)$ for some acyclic quiver $Q_E$ with components Dynkin quivers and  $\mc T_2^E \cong D^b(\coh(\XX(\rp'_E)))$ for some  $\rp'_E$ of Dynkin type  with $\rp'_E\prec \rp$   and $\Hom(\mc T_1^E, \mc T_2^E)=\Hom(\mc T_2^E, \mc T_1^E)=0$.  Thus we see that for each $\mc A \in C^{\{\rm Id\}}_{D^b(\coh(\XX(\rp')))}(D^b(\coh(\XX(\rp))))$ we have $\mc A \subset \langle \mc T_1^{E}, \mc T_2^{E}\rangle $ for some $E \in \mc E$ with the described properties of $\langle \mc T_1^{E}, \mc T_2^{E}\rangle $, furthermore we claim that $\mc A \subset   \mc T_2^{E} \cong   D^b(\coh(\XX(\rp'_E)))$, in particular $\mc A \in C^{\{\rm Id\}}_{D^b(\coh(\XX(\rp')))}(\mc T_2^E)$. Indeed, now $\mc A \cong \mc T(\rp')$ and  $\mc T(\rp')$ is the derived category of representations of a connected acyclic quiver, therefore  we have a full  exceptional collection $A_1,A_2,\dots, A_n$ in $\mc A $ such that for any $1\leq j <n$ $\langle A_1,\dots, A_j \rangle $ is not orthogonal to $\langle A_{j+1},\dots, A_n \rangle$, and since  $\langle \mc T_1^{E}, \mc T_2^{E}\rangle $ is an orthogonal decomposition  we have either $\{A_i\}_{i=1}^n \subset  \mc T_1^{E}$ or $\{A_i\}_{i=1}^n \subset  \mc T_2^{E}$, hence either $\mc A \subset  \mc T_1^{E}$ or $\mc A  \subset  \mc T_2^{E}$, finally since we have  $\KK$-linear equivalence $\mc T_1^E \cong D^b(Q_E)$ for some acyclic quiver $Q_E$ with components Dynkin quivers,  we can exclude the case  $\mc A \subset  \mc T_1^{E}$ by the fact that there are infinitely many indecomposable objects in $\mc A$ up to shifts and using \eqref{nice properties2} again. 
				
				In the case  ${\rm rank}\left (K_0(D^b(\coh(\XX(\rp'))))\right)<{\rm rank}\left (K_0(D^b(\coh(\XX(\rp))))\right)$ we deduce that there is a finite set of triangulated subcategories $\{ \mc T_2^E : E \in \mc E \}$ of $D^b(\coh(\XX(\rp)))$ such that  $\mc T_2^E \cong  D^b(\coh(\XX(\rp'_E))) $ with $\rp'_E$ of Dynkin type  and $\rp'_E\prec \rp$   for each $E \in  \mc E$    and such that  \begin{gather} \label{inclusion in finiteness corollary} C^{\{\rm Id\}}_{D^b(\coh(\XX(\rp')))}(D^b(\coh(\XX(\rp)))) \subset \bigcup_{E \in \mc E} C^{\{\rm Id\}}_{D^b(\coh(\XX(\rp')))}(\mc T_2^E). \end{gather} If $C^{\{\rm Id\}}_{D^b(\coh(\XX(\rp')))}(D^b(\coh(\XX(\rp))))\neq \emptyset$, then $C^{\{\rm Id\}}_{D^b(\coh(\XX(\rp')))}( D^b(\coh(\XX(\rp'_E))))\neq \emptyset $ for some $E \in \mc E$ and   by \eqref{help for finiteness 0} and   the induction assumption $\rp'\preccurlyeq  \rp'_E\prec \rp$. Furthermore again by the induction assumption  $ C^{\{\rm Id\}}_{D^b(\coh(\XX(\rp')))}(\mc T_2^E)$ is finite for each $E \in \mc E$ and therefore by \eqref{inclusion in finiteness corollary}  $C^{\{\rm Id\}}_{D^b(\coh(\XX(\rp')))}( D^b(\coh(\XX(\rp))))$ is finite.
	\epr
	\section{Remarks on  $\abs{C^{\Gamma}_{1}(D^b(\coh(\XX(\rp))))}$}  \label{remarks on the number}
Combining results in \cite{GL1} and \cite{LM}  we will study further   $\abs{C^{\Gamma}_{\mc T(\rp')}(\mc T(\rp))}$ for the case when $\rp'$  is the unique sequence of Dynkin type with  $\abs{\rp'}=3$. This is the case $\rp'=(1,1,1)$ and then $\mc T(\rp')=D^b(K(2))$.  With the terminology explained in the introduction  the set $C^{\{\rm Id\}}_{\mc T(\rp')}(\mc T(\rp))$ is the set of   nc  curves in $\mc T(\rp)$ of non-commutative  genus 1, i. e. this is  $C^{\{\rm Id\}}_{1}(\mc T(\rp))$. 	Due to Fact \ref{equivalences} it is enough to determine  $\abs{C^{\Gamma}_{1}(D^b(\coh(\XX(\rp))))}$.

\begin{remark} \label{remark for exceptianal pairs} For any $\mc B \in C^{\{\rm Id\}}_{1}(D^b(\coh(\XX(\rp))))$ there  exists an exact  equivalences $D^b(\coh(\PP^1))\cong D^b(K(2))\cong \mc B$ and from \cite[(25), (26), (27)]{DK41} we see that there is a sequence $\{s_i\}_{i\in \ZZ}$ of exceptional objects in $\mc B$ such that 
\begin{gather} \label{list of exceptional pairs} \forall i \in \ZZ \qquad \qquad  \mc B = \langle s_i,s_{i+1} \rangle,  \ \   \hom(s_i,s_{i+1})=2,\  \hom^k(s_i,s_{i+1})=0 \ \mbox{for} \ k\neq 0, \\ \label{complete lists  of  exceptional} 
\textit{complete lists  of  exceptional  pairs and objects (up to shifts) are} \  \{(s_i,s_{i+1})\}_{i\in\ZZ} \ \textit{and} \ \{s_i\}_{i\in\ZZ}. \\
\label{triangles} \textit{there is an exact triangle} \ \ \  s_{i-1}\rightarrow  s_{i}^2 \rightarrow s_{i+1} \rightarrow  s_{i-1}[1] \ \ \ \textit{for any} \ \ i\in \ZZ.  \end{gather}
\end{remark}

In  \cite{LM} for any exceptional pair $(A,B) \in \coh(\XX(\rp))$  is studied the  smallest exact extension-closed subcategory of  $\coh(\XX(\rp))$ containing $A$ and $B$. We will denote this subcategory by $\mc Z (X,Y)$. We note first:

\begin{lemma} \label{exceptionally generated subcat}  Let $\mc A \subset \coh(\XX(\rp))$ be an  exact and extension-closed subcategory such that there exists a $\KK$-linear exact equivalence $\coh(\PP^1)\rightarrow \mc A$. Then:
	
	{\rm (a)} object $X$/ (pair of objects  $X,Y$)   in $\mc A$ is exceptional object/(exceptional pair) with respect to $\mc A$ iff it is exceptional object/(exceptional pair) with respect to $\coh(\XX(\rp))$.
	
{\rm (b)} For any exceptional pair $X,Y \in \mc A$ we have $\hom(X,Y)=2$, $\hom^k(X,Y)=0$ for $k\neq 0$, and  $\mc Z(X,Y)=\mc A$.

{\rm (c)} There is a sequence $\{s_i\}_{i\in \ZZ}$ of exceptional objects contained in  $\mc A$, such that each exceptional object in  $\mc A$ is isomorphic to  $s_i$ for some $i\in \ZZ$ and  each  exceptional pair $(X,Y)$ in $\mc A$ is isomorphic   to $(s_i, s_{i+1})$ for some $i\in \ZZ$. 
Furthermore we have equality  $\langle \{s_i\}_{i\in \ZZ} \rangle =\langle s_j, s_{j+1} \rangle $ in $D^b(\coh(\XX(\rp))$  for any $j\in \ZZ$.

{\rm (d)} From {\rm(b), (c)} follows that   $\langle \mc A_{exc} \rangle  \in C^{\{\rm Id\}}_{1}(D^b(\coh(\XX(\rp))))$, where  $\mc A_{exc}$ is the set of exceptional objects in $\mc A$.  
\end{lemma}
\bpr Since $\mc A$ is exact and extension-closed we have \eqref{Homs and exts the same} and hence (a) follows.

(b)  Let  $X,Y \in \mc A$ be an exceptional pair.
Since $\mc A$ is exact and extension-closed,  $\mc Z(X,Y)\subset \mc A$.

 Let $F$ be an equivalence functor  $\coh(\PP^1)\rightarrow \mc A$.  We will show that $\mc Z(X,Y) \supset \mc A$ by showing that for each $A\in {\rm vect}(\PP^1)$ and each $A\in \coh_0(\PP^1)$ holds $F(A) \in \mc Z(X,Y)$ and using Fact \ref{direct sum tor and torfree}.

   All exceptional pairs in $\coh(\PP^1)$ are of the form $(\mc O(n) , \mc O(n+1))$ (see e.g. \cite[the proof of Theorem 3]{LM}) hence   we have $X \cong F(\mc O(n))$, $Y\cong  F(\mc O(n+1))$ for some $n\in \ZZ$, in particular $ F(\mc O(n)),   F(\mc O(n+1)) \in \mc Z(X,Y)$ and  $\hom(X,Y)=2$, $\hom^k(X,Y)=0$ for $k\neq 0$.  The sequence $\{s_j: j\in \ZZ\}$ in $D^b(\coh(\PP^1))$  from Remark \ref{remark for exceptianal pairs} can be chosen to be $\{\mc O(j): j \in \ZZ \}$ hence we have a short exact sequence $ 0 \rightarrow \mc O(i-1)\rightarrow  \mc O(i)^2 \rightarrow  \mc O(i+1) \rightarrow  0 $ for each $i$. Using these short exact sequences, the fact that the  functor $F$ is exact, and that $\mc Z(X,Y)$ is exact an extension closed by induction one shows that  $F(\mc O(j))\in \mc Z(X,Y)$ for each $j\in \ZZ$.  From \cite[ Proposition 2.6]{GL} it follows that for each vector bundle $A\in {\rm vect}(\PP^1)$ we have $F(A) \in \mc Z(X,Y)$.

 To show that for  each $A\in \coh_0(\PP^1)$ holds $F(A) \in \mc Z(X,Y)$ it is enough to show that for  each simple sheaf $S_x \in \coh_0(\XX(\rp))$, $x\in \PP^1$ holds  $F(S_x) \in \mc Z(X,Y)$ because of  \cite[ Proposition 2.5]{GL}.  Each simple sheaf $S_x \in \coh_0(\XX(\rp))$ for each $x\in \PP^1$  arises as the co-kernel term of an exact sequence  $ 0 \rightarrow \mc O\rightarrow  \mc O(1) \rightarrow  S_x \rightarrow  0 $ (see again in the  \cite[Theorem 3]{LM}). Therefore by the already proved part and  the fact that the  functor $F$ is exact, and that $\mc Z(X,Y)$ is exact an extension closed  we get $F(S_x)\in \mc Z(X,Y)$ for each $x\in \PP^1$. So we proved (b).
 
 {\rm (c)} The sequence $\{s_i=F(\mc O(i))\}_{i\in \ZZ}$ satisfies the desired properties, since $F$ is equivalence and using (a).  The short exact sequence  $ 0 \rightarrow \mc O(i-1)\rightarrow  \mc O(i)^2 \rightarrow  \mc O(i+1) \rightarrow  0 $ give rise to a distinguished triangle  $ s_{i-1}\rightarrow  s_{i}^2 \rightarrow s_{i+1} \rightarrow  s_{i-1}[1]$ in $D^b(\coh(\XX(\rp)))$ for each $i\in \ZZ$ and it follows that  $ \langle s_{i},s_{i+1}\rangle = \langle s_{i-1},s_{i}\rangle  $ for each $i$, therefore  $\langle \{s_i\}_{i\in \ZZ} \rangle =\langle s_j, s_{j+1} \rangle $ for any $j\in \ZZ$.
\epr

\begin{lemma} \label{from cohP1 tp genus 1 ncc} Let $\rp$ be of Dynkin type. 
 For any $\mc B \in C^{\{\rm Id\}}_{1}(D^b(\coh(\XX(\rp))))$ there  exists  a unique  exact and extension-closed subcategory  $\mc A \subset \coh(\XX(\rp))$,  such that  there exists a $\KK$-linear exact equivalence $\coh(\PP^1)\rightarrow \mc A$ and $\mc B = \langle \mc A_{exc} \rangle$, where  $\mc A_{exc}$ is the set of exceptional objects in $\mc A$.
\end{lemma}
\bpr  Let $\mc B \in C^{\{\rm Id\}}_{1}(D^b(\coh(\XX(\rp))))$ and let $\{s_i\}_{i\in \ZZ}$ be the sequence from Remark \ref{remark for exceptianal pairs}.
Using that $\coh(\XX(\rp))$ is hereditary we get integer $k_i\in \ZZ$, s.t.  $s_i[k_i]\in \coh(\XX(\rp))$ for any $i\in \ZZ$  and \eqref{list of exceptional pairs} shows that either $k_i=k_{i+1}$ or $k_{i}=k_{i+1}+1$. We will show that must hold $k_i=k_{i+1}$. Take any $i\in \ZZ$ and denote for short $A_1=s_i[k_i]$, $A_2=s_{i+1}[k_{i+1}]$. Let $\mc Z(A_1,A_2)$ be  the smallest  exact  and extension-closed  subcategory in $\coh(\XX(\rp))$ containing $A_1,A_2$.  
We claim that $\mc Z(A_1,A_2)$ is equivalent to $\coh(\PP^1)$ and that $\hom(A_1,A_2)=2$,  and then from \eqref{list of exceptional pairs} it follows that $k_i=k_{i+1}$.

We  show first   that $\mc Z(A_1,A_2)$ is not equivalent to a module category over a finite dimensional hereditary algebra. 
Indeed, if $\mc Z(A_1,A_2)$  is equivalent to a module category over a finite dimensional  hereditary algebra,  from the arguments as those related to \eqref{Homs and exts the same}, \eqref{Hom 2 or ext 2} follows that   $\mc Z(A_1,A_2) $ is equivalent to  $Rep_\KK(K(r))$ for some $r\geq 0$.  Since all exceptional pairs in $Rep_\KK(K(r))$ have $r$-dimensional Hom or Ext between their elements, using that either $\hom(A_1, A_2)=2$ or $\hom^1(A_1,A_2)=2$ due to \eqref{list of exceptional pairs} we conclude that $r=2$, however this contradicts \cite[Proposition 4.3]{LM}. 

Now,   from \cite[Theorem 2]{LM} we get   that  $\mc Z(A_1,A_2)$  is equivalent to $\coh(\PP^1)$, therefore $(A_1, A_2)$ corresponds to an exceptional pair in  $\coh(\PP^1)$, however all exceptional pairs in $\coh(\PP^1)$ are of the form $(\mc O(n) , \mc O(n+1))$ (see e.g. \cite[the proof of Theorem 3.]{LM}), and therefore $0\neq \hom(A_1,A_2)=2$.   

We see that $k_i=k_{i+1}$ for each $i\in \ZZ$ and it follows that we can ensure 
\begin{gather} \label{for all i si in coh}
\forall i \in \ZZ \ \ \ s_i\in \coh(\XX(\rp)).
\end{gather}
From \eqref{triangles},   it follows that  for any two integers $i,j\in \ZZ$ we have
\begin{gather} \label{help for uniqueness}
\mc Z(s_i,s_{i+1}) = \mc Z(s_j,s_{j+1}). 
\end{gather}
If we denote $\mc A= \mc Z(s_i,s_{i+1}) $, then \eqref{help for uniqueness}  implies that $ s_i\in \mc A$ for any  $i \in \ZZ$. Recalling that  $(s_i,s_{i+1})$ is an exceptional pair for each $i$, seeing  the properties of the sequence  from Lemma \ref{exceptionally generated subcat} (c), and   taking into account  Fact \ref{make an element unique} we see that in fact $\{s_i\}_{i\in \ZZ}$ is a sequence as in (c)  Lemma \ref{exceptionally generated subcat}, hence  the set of all exceptional objects in $\mc A$  up to isomorphism is $\{s_i\}_{i\in \ZZ}$ and $\mc B = \langle  \mc   A_{exc} \rangle=\langle s_i, s_{i+1} \rangle$.

Now we prove the uniqueness.  Assume that  $\mc B \in C^{\{\rm Id\}}_{1}(D^b(\coh(\XX(\rp))))$ and that  $\mc B= \langle \mc  A_{1exc}  \rangle = \langle \mc  A_{2exc}  \rangle$  where  $\mc A_i$   is an exact,  extension closed subcategory  $\mc A_i \subset \coh(\XX(\rp))$,  such that  there exists a $\KK$-linear exact equivalence $\coh(\PP^1)\rightarrow \mc A_i$ for $i=1,2$. From Lemma \ref{exceptionally generated subcat} (c) we see that     $\mc B= \langle  X_1,Y_1  \rangle = \langle X_2,Y_2 \rangle$, where   $(X_i,Y_i)$ is an exceptional pair in $\mc A_i$ for $i=1,2$.

   From Lemma \ref{exceptionally generated subcat} (b) it follows that for $i=1,2$ we have   $\mc A_i = \mc Z(X_i,Y_i)$.  The pair $(X_1,Y_1)$ is  exceptional  in $\mc B$, whose elements are in $\coh(\XX(\rp))$, therefore   \eqref{complete lists  of  exceptional} and \eqref{for all i si in coh} imply that  it is  isomorphic to a pair of the form $(s_i,s_{i+1})$ for some $i\in \ZZ$, similarly $(X_2,Y_2)$ is isomorphic to a pair of the form $(s_j,s_{j+1})$ for some $j\in \ZZ$. Now from \eqref{help for uniqueness} we get    $\mc A_1 =\mc A_2$. The lemma is proved.
\epr

	\begin{coro} \label{with the full group} Let  $\rp=(p_1,p_2,p_3)$ be any  weight sequences of Dynkin type as defined in Definition  \ref{the categories indexed by p}. 
				Then  $\abs{C^{{\rm Aut(\mc  T(\rp))}}_{1}(\mc T(\rp))}= 1$.
				
					It follows that for any  affine connected acyclic quivers  $Q$ holds $\abs{C^{{\rm Aut(D^b(Q))}}_{1}(D^b(Q))}=1$.
				\end{coro}
\bpr 	The last sentence follows from Fact \ref{equivalences} and the comments in Definition \ref{the categories indexed by p}. 

Let $\mc B \in C^{\{\rm Id\}}_{1}(D^b(\coh(\XX(\rp))))$. To prove the corollary  we will find  $\beta \in{\rm Aut}(D^b(\coh(\XX(\rp)))) $ such that $\beta(\mc B)=\langle \mc Z(\mc O,\mc O(\textbf{c}))_{exc} \rangle $, where $(\mc O,\mc O(\textbf{c}))$ is the corresponding  subsequence of the sequence in  Fact \ref{canonical exceptional collection}. 
 From Lemma \ref{from cohP1 tp genus 1 ncc} we have $\mc B=\langle \mc A_{exc} \rangle $ for some exact extension closed subcategory $\mc A$ of  $ \coh(\XX(\rp))$ which is equivalent to $\coh(\PP^1)$.  From Lemma \ref{exceptionally generated subcat} (b) we deduce that  $\mc  A= \mc Z(X,Y)$ for some exceptional pair $(X,Y)$. Since $\mc Z(X,Y)$ is not equivalent to a module category ( see also \cite[Corollary 1]{LM}),  from \cite[Theorem 3]{LM} it follows that for some auto-equivalence $\alpha \in {\rm Aut}(\coh(\XX(\rp)))$ holds $\alpha \left (\mc Z(X,Y)\right )= \mc Z(\mc O,\mc O(\textbf{c}))$, this auto-equivalence gives rise to an auto-equivalence $\beta \in{\rm Aut}(D^b(\coh(\XX(\rp)))) $ s.t.  $\mc Z(\mc O,\mc O(\textbf{c}))_{exc}$ is the isomorphism closure of the set $\beta \left (\mc A_{exc}\right )$,  and it follows that $\beta(\mc B)=\langle \mc Z(\mc O,\mc O(\textbf{c}))_{exc} \rangle  $.
\epr 
	In  \cite[COROLLARY 9.9]{GL1}  is written that ``\textit{up to equivalence of functors there are exactly $p_0 \cdots p_n$ exact embeddings $\Phi_t$, $t=0,\dots, p_0 \cdots p_n$ from  $\coh(\PP^1)$ to $\coh(\XX(\rp))$, whose image is closed under extensions.} '' We could not find in \cite[COROLLARY 9.9]{GL1} an explanation,  what is the equivalence  meant by saying ``\textit{up to equivalence of functors}''. In a previous version of our paper   we  presumed  that
	
	 \begin{gather} \mbox{\textit{two exact  embeddings $F,G$ from $\coh(\PP^1)$ to $\coh(\XX(\rp))$  should be considered equivalent}}  \nonumber \\ \label{hypothesis}
	 \mbox{ \textit{in \cite[COROLLARY 9.9]{GL1}  iff $F\circ \alpha$ is naturally isomorphic to $G$ }} \\ \mbox{\textit{for some $\KK$-linear auto-equivalence $\alpha$ of $\coh(\PP^1)$.} } \nonumber
	  \end{gather}
	  
	  Professor Helmut Lenzing confirmed in a private communication  that \eqref{hypothesis} is correct and now we can infer: 
	
	\begin{remark} \label{Remark on p_1 dots p_t} Let $\rp =(p_1, p_2, p_3)$ be of Dynkin type  as defined in Definition  \ref{the categories indexed by p}..  Since \eqref{hypothesis} holds then by similar arguments as in the proof of  \cite[Lemma 4.5]{DK5} one can show that \cite[COROLLARY 9.9]{GL1} implies that the set of    exact, extension-closed subcategories in $\coh(\XX(\rp))$, which are $\KK$-equivalent to $\coh(\PP^1)$ contains  $p_1 p_2 p_3$ elements.
		
	\end{remark} 
\begin{coro}  \label{coro for numbers} For any sequence $\rp =(p_1, p_2, p_3)$ of Dynkin type (as defined in Definition \ref{the categories indexed by p}) holds $\abs{C^{\{\rm Id\}}_{1}(\mc T(\rp))}= p_1 p_2 p_3$.
	
\end{coro}
 \bpr Follows from   Lemmas \ref{from cohP1 tp genus 1 ncc},  \ref{exceptionally generated subcat} (d)  and Remark \ref{Remark on p_1 dots p_t}.  \epr
 
 \section{Final remarks and conjectures} \label{final}
 We expect  that with the help of the cited here   works of  Geigle, Lenzning, Meltzer, Hübner   finiteness  of  other  sets of the type $C^{\{\rm Id \}}_{D^b(\coh(\XX'))}(D^b(\coh(\XX)))$ can be proved, where $\XX'$ and  $\XX$ are weighted projective lines.

Another  example studied   in \cite{DK5} is $D^b(\PP^2)$, here we fix $\KK=\CC$. In particular it is shown that then $C_{-1}^{\{\rm Id \}}(\mc T) =\emptyset$, and  $\forall l\geq 0$ the set $C_{l}^{\langle S \rangle}(\mc T)$ is finite, and:\footnote{Recall that a Markov number $x$ is a number $x \in \NN_{\geq 1}$ such that there exist integers $y,z$ with $x^2+y^2+z^2=3 x y z$.} 
\begin{gather} \label{non-rivials for markov} \left \{l\geq 0:C_{l}^{\langle S \rangle}(\mc T) \neq \emptyset \right \}= \left \{l\geq 0:C_{l}^{\{\rm Id \}}(\mc T) \neq \emptyset \right \}= \{3 m - 1: m \ \mbox{is a Markov number} \}.
\end{gather}

Thus,   there are infinitely many non-commutative curves in $D^b(\PP^2)$. There are also infinitely many  smooth projective curves in $\PP^2$, in contrast to that we have:

\begin{lemma} \label{no commuta in P2} For any smooth projective curve $S$   we have $C_{D^b(S)}^{\{\rm Id\} }(D^b(\PP^2))=\emptyset$. 
	
\end{lemma} \bpr We will omit writing  the superscript $\{\rm Id\}$ in $C_{D^b(S)}^{\{\rm Id\} }(D^b(\PP^2))$. 
From \eqref{non-rivials for markov} it follows that $C_{1}(D^b(\PP^2))=\emptyset$  and hence  $C_{D^b(\PP^1)}(D^b(\PP^2))=\emptyset$, since $D^b(\PP^1)\cong D^b(K(2))$ by \cite{B}.

 Let  $S$ be a smooth projective curve of positive   genus. We will show that there is no a full triangulated subcategory in $D^b(\PP^2)$ equivalent to $D^b(S)$ and then from Lemma \ref{bijection in formula} it follows that $C_{D^b(S)}^{\{\rm Id\} }(D^b(\PP^2))=\emptyset$. Due to the fact that $D^b(\PP^2)$ is generated by an exceptional triple we have 
 $K_0(D^b(\PP^2))\cong \ZZ^3$ and hence:
 \begin{gather} \label{help fact} \mbox{\textit{each subgroup of $K_0(D^b(\PP^2))\cong \ZZ^3$ is free of rank at most $3$.}} \end{gather}
 
 If there is a full triangulated subcategory  $\mc B \subset D^b(\PP^2)$ which is equivalent to $D^b(S)$ then there is a SOD $\langle \mc B^{\perp} , \mc B \rangle$ of $D^b(\PP^2)$ (see \cite[page 3]{BO}) and then we would have that $K_0(\mc B)\cong K_0(D^b(S))$ is a subgroup of $K_0(D^b(\PP^2))$. Therefore from \eqref{help fact} it follows that $ K_0(D^b(S))$ is free of finite rank.   However the Grothendieck group of a smooth projective curve $S$ is isomorphic to $Pic(S)\oplus \ZZ$ (see \cite[Remark on page 34]{Potier}) and $Pic(S)$ contains a subgroup, $Pic_0(S)$, isomorphic to the Jacobian variety of the curve $S$ (see e.g. \cite[p. 168]{Waldschmidt}) and for positive genus this is not a group of finite rank (it is even uncountable). 
\epr

Finally we discuss cases of projective threefolds for which there are smooth curves $S$, such that $C_{D^b(S)}(D^b(X))\neq \emptyset$.

\begin{conj} \label{final conjecture}
 From \cite[Corollary 3.9, Remarks 3.8, 3.12]{BerBolo} it follws that  for a threefold $X$
 satisfying certain properties, in particular its intermediate Jadobian splits in Jacobians of curves $ J(X)=\bigoplus_{i=1}^N J(S_i)$,
  and for a smooth projective curve $S$ of positive genus  there exists a fully faithful functor  $D^b(S)\rightarrow D^b(X)$ iff $S\cong S_i$ for some $i\in \{1,\dots, N\}$.  Therefore for such threefold $X$ we have $C_{D^b(S)}^{\{\rm Id\} }(D^b(X))\neq \emptyset$ iff  $S\cong S_i$ for some $i\in \{1,\dots, N\}$. Examples of $X$ for which  this holds are (see \cite[ Remark 3.12]{BerBolo}): 
  
  (a) Smooth quadric,   or a $\PP^1$-bundle over a rational surface, or a $\PP^2$ -bundle over $\PP^1$ , or a $V_5$, or a  $V_2$ Fano threefold.
  
   (b) $X$ is   the complete intersection of
   two quadrics or a Fano threefold of type $V_{18}$, and  $J(\Gamma) \cong  J(X)$  with $\Gamma$ hyperelliptic
   
   (c) Rational conic bundle $X\rightarrow S$     over a minimal surface.
   
   (d) Del Pezzo fibrations: $X \rightarrow \PP^1$ is a quadric fibration with at most simple degenerations,
   in which case the hyperelliptic curve $\Gamma \rightarrow \PP^1$ ramified along the degeneration appears
   naturally as the orthogonal complement of an exceptional sequence of $D^b (X)$.
   $X \rightarrow \PP^1$ is a rational Del Pezzo fibration of degree four.

   We conjecture that  $\abs{C_{D^b(S_i)}^{{\rm Aut}(D^b(X)) }(D^b(X))}<\infty$  for such $X$ and any $i=1,\dots, N$.  
 \end{conj} 
\let\oldaddcontentsline\addcontentsline
\renewcommand{\addcontentsline}[3]{}

\let\addcontentsline\oldaddcontentsline

\end{document}